\documentclass[twoside,12pt]{article}
\usepackage{t1enc}
\usepackage[latin2]{inputenc}

\usepackage{amssymb,amsmath}

\newenvironment{proof}{\begin{trivlist}\item[]{\it
Proof.}}{\hfill$\square$\end{trivlist}}

\newtheorem{theorem}{Theorem}[section]
\newtheorem{corollary}[theorem]{Corollary}

\newcommand{\C}{\mathbb C}
\newcommand{\la}{\lambda}

\newcommand{\tr}{\mathrm {tr}\;}

\newcommand{\diag}{\mathrm {diag}}
\newcommand{\Hom}{\mathrm {Hom}}
\newcommand{\Z}{\mathbb Z}

\begin{document}

\date{}

\title{ Character formulae for classical groups }

\author{ P. E. Frenkel
\thanks{ Partially supported by OTKA grants T 042769 and T 046365.}}

\maketitle 
\begin{abstract} 
We give formulae relating the value $\chi_\la \left(g\right)$ of an irreducible character
of a classical group $G$ to  entries of  powers of the matrix $g\in G$.
This 
yields a  far-reaching generalization   of a 
result of J.\ L.\ Cisneros-Molina  concerning the $GL_2$ case [C]. 
\end{abstract} 



\bigskip

\section*{Introduction}
The Weyl character formula tells us how to compute the character $\chi_\lambda$
of an irreducible finite dimensional representation $V_\lambda$ with highest weight
$\lambda$
of a (complex, 
semisimple, connected) Lie group $G$:  $$\sum_{\sigma\in\mathfrak W} \left(-1\right)^\sigma
  z^{\sigma\left(\lambda+\rho\right)}=\chi_\lambda\cdot 
 \sum_{\sigma\in\mathfrak W} \left(-1\right)^\sigma
  z^{\sigma\rho}.$$
Here, for each weight $\ell\in \mathfrak h^*$ of the Lie algebra $\mathfrak g$
of $G$, the exponential   $z^\ell :\tilde H\to\C^*$ is the corresponding
multiplicative character of the preimage $\tilde H$ of a maximal torus $H\le
G$
in a universal covering
$\tilde G\to G$.
The weight $\rho\in\mathfrak h^*$ is the half-sum of the positive roots, $\mathfrak W$ is the Weyl
group. Note that 
both sides of the formula are $\mathfrak W$--antisymmetric characters of
$\tilde H$, but $\chi_\lambda$ is well-defined as a $\mathfrak W$--symmetric
character of $H$.  Note also that $\Delta=\sum_{\sigma\in\mathfrak W} \left(-1\right)^\sigma
  z^{\sigma\rho}$ is not identically zero, so $\chi_\lambda$
 is  expressed as a ratio of Laurent polynomials in the coordinates on
$\tilde H$. (We know {\it a priori} that $\chi_\lambda$ is itself a  Laurent
 polynomial in the coordinates on $H$, but
  with many more terms in general than the numerator  and denominator.)

The formula expresses the value of the character $\chi_\lambda$
at a group element $g\in G$ in terms of 
a conjugate of the semisimple part of  $g$ 
in  the  maximal torus $H$, i.e., in the case of the classical matrix groups,
in terms of
 the eigenvalues of $g$. There are  equally explicit expressions, called
determinantal identities or Giambelli formulae [FH, Section A.1, formulae
  (A.5), (A.6) and Section A.3], in terms of  the
elementary resp.\ the complete symmetric polynomials in the eigenvalues.
 
In the present paper, we
consider the connected classical groups and we prove 
variants of the Weyl character formula  that express
the value $\chi_\lambda\left(g\right)$  in many different, explicit
 rational ways in terms of
 the entries of powers of the generic matrix $g\in G$. 
It seems likely that these
 formulae provide the fastest and most straightforward
way of calculating $\chi_\la (g)$ for generic $g$.


This paper was motivated by J.\
L.\ Cisneros-Molina's paper [C]  whose main result is the following. Let
$\omega\ne 0$
be  a
linear function on the space $M_2
$ of 
$2\times 2$ matrices, such that $\omega\left(\mathbf 1\right)=0$. For example, $\omega$ could be one of the
two off-diagonal entries.  Then, for $\lambda=0,1,\dots$, we have $\omega\left(g^{\lambda+1}\right)/\omega\left(g\right)=\tr S^\lambda g$, the
trace of the action of $g$ on the  $\lambda$-th symmetric power of the
standard
vector
representation. 
Our results can be
considered as far-reaching generalizations of this fact.
In particular,  for $g\in M_{r+1}$, and $\la=0,1,\dots$, the trace $\tr S^\la
g$ equals the ratio of the $r$--dimensional
volumes of the two parallelepipeds spanned in
$M_{r+1}
/M_1
$ 
by
the images of  $g$, $g^2$, \dots, $g^{r-1}$,
$g^{\la+r}$ and of  $g$, $g^2$, \dots, $g^{r-1}$,
$g^{r}$, respectively 
(except when both volumes are zero, i.e., $g$ has a minimal polynomial of
degree $<r+1$). This is a particular  case of Corollary~\ref{o} below.

Our proofs are motivated by the first of the four proofs given in [C], which
is  due to Jeremy Rickard.


\section{General linear group}
  Let $G=GL_{r+1}\left(\C\right)$ with $H=
  \left(\C^*\right)^{r+1}$ the  maximal torus consisting of all invertible
diagonal matrices,
  $
\Hom \left(H,\C^*\right)= \Z^{r+1}$ the weight lattice,
and $\mathfrak W=\mathfrak S_{r+1}$
  the Weyl group.  Write $$\rho=\left(\frac r2,\frac{r-2}2,\dots,\frac{2-r}2,
  \frac {-r}2\right)$$ for the
  half-sum of the positive roots. Set $\rho_t=\rho +\left(t,t,\dots,t,t\right)
$ for $t\in\C$.
For 
$\lambda=\left(\lambda_0,\dots,\lambda_l\right)\in\Z^{r+1}$, write $z^\lambda :H\to\C^*$
for the corresponding multiplicative character of the torus $H$, and, when
$\lambda$ is dominant, i.e. $\lambda_0\ge\dots\ge\lambda_r$, write $\chi_\lambda:G\to\C$
for the character of the irreducible representation with highest weight
$\lambda$.
The Weyl character formula   $$\sum_{\sigma\in\mathfrak W} \left(-1\right)^\sigma
  z^{\sigma\left(\lambda+\rho_t\right)}=\chi_\lambda 
 \cdot\sum_{\sigma\in\mathfrak W} \left(-1\right)^\sigma
 z^{\sigma\rho_t}$$ holds with any $t\in\C$.
The freedom in the choice of $t$ comes from the central $\C^*$ in $G$.
Both sides of the formula 
are $\mathfrak W$--antisymmetric characters of the infinite cover $\tilde H$. When
$\rho_t\in\Z^{r+1}$, both sides descend to $H$. 

For $\ell\in\C^{r+1}$ and $g\in G$, 
define
  $$g^\ell=\bigotimes_{i=0}^rg^{\ell_i}\in M_{r+1}\left(\C\right)^{\otimes \left(r+1\right)}.$$
This is multi-valued, it
depends on a choice of the value of $\log g\in \mathfrak{gl}_{r+1}\left(\C\right)$.  When
$\ell\in\Z^{r+1}$, it is  single-valued.
When $\ell\in\Z_{\ge 0}^{r+1}$, we may allow  $g\in M_{r+1}\left(\C\right)$ rather
than $g\in G$.

We have $$\sum_{\sigma\in\mathfrak W} \left(-1\right)^\sigma
  g^{\sigma\ell}=\bigwedge_{i=0}^{r} g^{\ell_i}\in M_{r+1}\left(\C\right)^{\wedge\left(r+1\right)}.$$


\begin{theorem}\label{GL}  
Let
  $\lambda=\left(\lambda_0\ge\dots\ge\lambda_r\right)\in\Z^{r+1}$  and $g\in
  GL_{r+1}\left(\C\right)$. Then, for any $t\in\C$,
 $$\sum_{\sigma\in\mathfrak W} \left(-1\right)^\sigma
  g^{\sigma\left(\lambda+\rho_t\right)}=\chi_\lambda \left(g\right)
 \cdot
\sum_{\sigma\in\mathfrak W} \left(-1\right)^\sigma
  g^{\sigma\rho_t};$$
equivalently,
$$\bigwedge_{i=0}^{r} g^{\lambda_i+r/2-i+t}=\chi_\lambda
  \left(g\right)\cdot\bigwedge_{i=0}^{r} g^{r/2-i+t},$$ where the powers are defined using any
  (but always the same) value of $\log g$.  When $\rho_t\in\Z^{r+1}$, the
  powers are single-valued. In particular, for $t=r/2$, we get$$\bigwedge_{i=0}^{r} g^{\lambda_i+r-i}=\chi_\lambda
  \left(g\right)\cdot\bigwedge_{i=0}^{r} g^{r-i}.$$  When $\rho_t$ and $\lambda$ are
  both in $\Z_{\ge 0}^{r+1}$, we may allow  $g\in M_{r+1}\left(\C\right)$.
\end{theorem}

\begin{proof}
The set of diagonalizable invertible matrices 
is dense in $M_{r+1}\left(\C\right)$, so we may
assume that $g$ is such.  The statement of the theorem is invariant under
conjugation, so we may assume that $g=\diag \left(z_0,\dots,z_r\right)\in H$. Then 
$$\sum_{\sigma\in\mathfrak W} \left(-1\right)^\sigma
  g^{\sigma\ell}=\bigwedge_{i=0}^{r}
  g^{\ell_i}=\left|z_j^{\ell_i}\right|
\cdot
\bigwedge_{j=0}^r e_{jj},$$
where $e_{jj}$ is the diagonal matrix with a single 1 at the $j$-th position.
The theorem now follows from the Weyl character formula.
\end{proof}



\begin{corollary}\label{O} 
Let $\Omega$ be an alternating $\left(r+1\right)$--linear form on the space
$M_{r+1}\left(\C\right)$. Then,
for $\lambda=\left(\lambda_0\ge\dots\ge\lambda_r\right)\in\Z^{r+1}$, we have 
 $${\Omega \left( g^{\lambda_0+r}, g^{\lambda_1+r-1},  
\dots, g^{\lambda_r}\right)}=\chi_\lambda \left(g\right)\cdot {\Omega
\left(g^r, g^{r-1}, \dots, \mathbf 1\right)}.$$
\end{corollary} 
To express  $\chi_\la(g)$ as  a rational function in  entries of powers of 
$g$, we must choose $\Omega$ such that the right hand side is not identically
zero. For example, $\Omega (g_0,\dots, g_r)$  could be the determinant of the
matrix formed by the diagonals, or by the first rows, etc.\
 of the argument matrices.  
To calculate $\chi_\la (g)$ 
for a numerically given $g$, we need to choose $\Omega$ such that
$\Omega(g^r,g^{r-1},\dots,{\bf{ 1}})\ne 0.$ This is possible if and only if
$\bigwedge_{i=0}^r g^{r-i}\ne 0$, i.e., the minimal polynomial of $g$ is  its
characteristic polynomial.
When
$g$ has a minimal polynomial of lower degree, we can use l'Hospital's rule.


\begin{corollary}\label{o}
Let $\omega$ be an alternating $r$--linear form on the space
$M_{r+1}\left(\C\right)$ such that $\omega$
vanishes if an argument is $\bf 1$.
Then, for $\lambda$ as above and with $\lambda_r=0$, we have
 $${\omega \left( g^{\lambda_0+r}, g^{\lambda_1+r-1},  \dots,
 g^{\lambda_{r-1}}\right)}=\chi_\lambda \left(g\right)\cdot {\omega \left(g^r, g^{r-1}, \dots,g\right)}.$$
\end{corollary}
To express  $\chi_\la(g)$ as  a rational function in  entries of powers of 
$g$, we must choose $\omega$ such that the right hand side is not identically
zero. For example, $\omega (g_0,\dots, g_{r-1})$  could be the determinant of the
$r\times r$
matrix formed by the truncated (i.e., leftmost entry omitted)
first rows
 of the argument matrices.  
To calculate $\chi_\la ( g)$ 
for a numerically given $g$, we need to choose $\omega$ such that
$\omega(g^r,g^{r-1},\dots,g)\ne 0.$ This is possible if and only if
the minimal polynomial of $g$ is  its
characteristic polynomial.

Corollary~\ref{o}, for $r=1$, is the result of  J.\ L.\ Cisneros-Molina's paper
  [C] mentioned in the Introduction.

To derive Corollary~\ref{o} from Corollary~\ref{O}, simply set $\Omega=d\omega$, defined
as usual by $$\Omega\left(g_0,\dots,g_{r}\right)=\sum_{i=0}^{r}\left(-1\right)^i\omega\left(g_0,\dots,
g_{i-1},g_{i+1},\dots, g_{r}\right).$$ Then $\Omega \left(g_0,\dots,g_{r-1},\mathbf 1\right)=\left(-1\right)^r\omega
 \left(g_0,\dots,g_{r-1}\right)$ and Corollary~\ref{o} follows.

\section{Special linear group}
  Let $G=SL_{r+1}\left(\C\right)$ with $H\simeq
  \left(\C^*\right)^{r}$ the  maximal torus consisting of all unimodular diagonal matrices,
  $
\Hom \left(H,\C^*\right)= \Z^{r+1}/\Z$ the weight lattice,
and $\mathfrak W=\mathfrak S_{r+1}$
  the Weyl group.  Write $\rho=\left(r,r-1,\dots,0\right)+\Z\cdot \left(1,\dots,1\right)\in\Z^{r+1}/\Z$
for  the half-sum of the positive  roots.
When 
$\lambda=\left(\lambda_0,\dots,\lambda_r\right)\in\Z^{r+1}/\Z$, write $z^\lambda :H\to\C^*$
for the corresponding multiplicative character of the torus $H$, and, when
$\lambda$ is dominant, write $\chi_\lambda:G\to\C$
for the character of the irreducible representation with highest weight
$\lambda$.
The Weyl character formula is valid as stated in the introduction.
Both sides are $\mathfrak W$--antisymmetric characters of $H$.

For $\ell\in\Z^{r+1}/\Z$ and $g\in G$, the antisymmetric tensor
  $$\bigwedge_{i=0}^{r} g^{\ell_i}\in M_{r+1}\left(\C\right)^{\wedge\left(r+1\right)}$$
is well defined because either $g$ has a minimal polynomial of degree $<r+1$,
in which case the algebra $\C[g]$ has dimension $<r+1$ and the antisymmetric
tensor above is zero, or else $g$ has its characteristic polynomial as minimal
polynomial, in which case $\dim \C[g]=r+1$ and multiplication by  $g$ on it has
determinant $\det g=1$, so the tensor is independent of the chosen
representative of $\ell$.

\begin{theorem}
Let
  $\lambda=\left(\lambda_0\ge\dots\ge\lambda_r\right)+\Z\cdot \left(1,\dots,1\right)\in\Z^{r+1}/\Z$  and $g\in SL_{r+1}\left(\C\right)$. Then
$$\bigwedge_{i=0}^{r} g^{\ell_i}=\chi_\lambda
\left(g\right)\cdot\bigwedge_{i=0}^{r} g^{r-i},$$
where $\ell_i=\lambda_i+r-i$.
\end{theorem}

\begin{proof}  The theorem trivially follows from Theorem~\ref{GL}.
\end{proof}




\section{Odd special orthogonal group}
  Let $G=SO_{2r+1}\left(\C\right)$ be the connected group preserving the quadratic form
$$x_1y_1+\dots+x_ry_r+z^2.$$ 
We take the maximal torus $H=
  \left(\C^*\right)^{r}$ consisting of all special orthogonal  diagonal 
matrices $\diag \left(z_1,z_1^{-1},\dots,z_r,z_r^{-1},1\right)$.
In  the weight lattice  $
\Hom \left(H,\C^*\right)= \Z^{r}$, we take
  $\lambda=\left(\lambda_1,\dots,\lambda_r\right)$ to correspond to the monomial
  $z^\lambda=\prod z_j^{\lambda_j}$.
The Weyl group  $\mathfrak W$ is the semidirect product of $\mathfrak S_{r}$
and $ Z_2^r$.  Write $\rho=\left(r-1/2,r-3/2,\dots,3/2,1/2\right)$
for  the
half-sum of the positive  roots.
The Weyl character formula is valid as stated in the introduction.
Both sides are $\mathfrak W$--antisymmetric characters of the 
 double cover $\tilde H$.

For $\ell\in\left(\Z +\frac 12\right)^{r}$ and $g\in G$, define
  $$g^\ell=\bigotimes_{i=1}^rg^{\ell_i}\in M_{2r}\left(\C\right)^{\otimes r}.$$
This is multi-valued, it depends on a  choice of $\sqrt g$.
We have $$\sum_{\sigma\in\mathfrak W} \left(-1\right)^\sigma
  g^{\sigma\ell}=\bigwedge_{i=1}^{r}
  \left(g^{
\ell_i}-g^{-\ell_i}\right)\in M_{2r}\left(\C\right)^{\wedge r}.$$

\begin{theorem}  
Let
  $\lambda=\left(\lambda_1\ge\dots\ge\lambda_r\right)\in\Z_{\ge 0}^{r}$  and $g\in SO_{2r+1}\left(\C\right)$. Then
 $$\sum_{\sigma\in\mathfrak W} \left(-1\right)^\sigma
  g^{\sigma\ell}=\chi_\lambda \left(g\right)\cdot
 \sum_{\sigma\in\mathfrak W} \left(-1\right)^\sigma
  g^{\sigma\rho};$$
equivalently,
$$\bigwedge_{i=1}^{r} \left(g^{\ell_i}- g^{-\ell_i}\right)=\chi_\lambda
\left(g\right)\cdot\bigwedge_{i=1}^{r} 
\left(g^{r+1/2-i}-g^{- \left(r+1/2-i\right)}\right),$$ where $\ell=\lambda+\rho$,
i.e. $\ell_i=\lambda_i+r+1/2-i$, and the powers are defined using
any, but always the same value of $\sqrt g\in SO_{2r+1}\left(\C\right)$.
\end{theorem}

\begin{proof}
The set of diagonalizable matrices is dense in $G$, so we may
assume that $\sqrt g$ is such.  The statement of theorem is invariant under
conjugation, so we may 
assume that $$\sqrt g=\diag \left(z_1^{1/2},z_1^{-1/2},\dots,z_r^{1/2},z_r^{-1/2},1\right)\in H.$$ Then 
$$\sum_{\sigma\in\mathfrak W} \left(-1\right)^\sigma
  g^{\sigma\ell}=\bigwedge_{i=1}^{r}
  \left(g^{\ell_i}-g^{-\ell_i}\right)=\left|z_j^{\ell_i}-z_j^{-\ell_i}\right|\cdot
\bigwedge_{j=1}^r \left(e_{jj}-f_{jj}\right),$$
where $e_{jj}$ resp.\ $f_{jj}$
is the diagonal matrix with a single 1 at the position
  corresponding to the $x_j$ resp.\ $y_j$ coordinate.
The theorem now follows from the Weyl character formula.
\end{proof}



\section{Symplectic group}
  Let $G=Sp_{2r}\left(\C\right)$ be the group preserving the skew bilinear form
  $\sum_{i=1}^r\left(x'_iy''_i-y'_ix''_i\right)$ on $\C^{2r}$. We take the maximal torus $H=
  \left(\C^*\right)^{r}$ consisting of all symplectic  diagonal matrices $\diag \left(z_1,z_1^{-1},\dots,z_r,z_r^{-1}\right)$.
In  the weight lattice  $
\Hom \left(H,\C^*\right)= \Z^{r}$, we take
  $\lambda=\left(\lambda_1,\dots,\lambda_r\right)$ to correspond to the monomial
  $z^\lambda=\prod z_j^{\lambda_j}$.
The Weyl group  $\mathfrak W$ is the semidirect product of $\mathfrak S_{r}$
and $ Z_2^r$.  Write $\rho=\left(r,r-1,\dots,1\right)$
for  the
half-sum of the positive  roots.
The Weyl character formula is valid as stated in the introduction.
Both sides are $\mathfrak W$--antisymmetric characters of $H$.

For $\ell\in\Z^{r}$ and $g\in G$, define
  $$g^\ell=\bigotimes_{i=1}^rg^{\ell_i}\in M_{2r}\left(\C\right)^{\otimes r}.$$
Then $$\sum_{\sigma\in\mathfrak W} \left(-1\right)^\sigma
  g^{\sigma\ell}=\bigwedge_{i=1}^{r}
  \left(g^{
\ell_i}-g^{-\ell_i}\right)\in M_{2r}\left(\C\right)^{\wedge r}.$$

\begin{theorem}
Let
  $\lambda=\left(\lambda_1\ge\dots\ge\lambda_r\right)\in\Z_{\ge 0}^{r}$  and $g\in Sp_{2r}\left(\C\right)$. Then
 $$\sum_{\sigma\in\mathfrak W} \left(-1\right)^\sigma
  g^{\sigma\ell}=\chi_\lambda \left(g\right)\cdot
 \sum_{\sigma\in\mathfrak W} \left(-1\right)^\sigma
  g^{\sigma\rho};$$
equivalently,
$$\bigwedge_{i=1}^{r} \left(g^{\ell_i}- g^{-\ell_i}\right)=\chi_\lambda
\left(g\right)\cdot\bigwedge_{i=1}^{r} 
\left(g^{r+1-i}-g^{- \left(r+1-i\right)}\right),$$ where $\ell=\la+\rho$, i.e., $\ell_i=\lambda_i+r+1-i$.
\end{theorem}

\begin{proof}
The set of diagonalizable matrices is dense in $G$, so we may
assume that $g$ is such.  The statement of theorem is invariant under
conjugation, so we may 
assume that $g=\diag \left(z_1,z_1^{-1},\dots,z_r,z_r^{-1}\right)\in H$. Then 
$$\sum_{\sigma\in\mathfrak W} \left(-1\right)^\sigma
  g^{\sigma\ell}=\bigwedge_{i=1}^{r}\left(
  g^{\ell_i}-g^{-\ell_i}\right)=\left|z_j^{\ell_i}-z_j^{-\ell_i}\right|
\cdot\bigwedge_{j=1}^r \left(e_{jj}-f_{jj}\right),$$
where $e_{jj}$ resp.\ $f_{jj}$
is the diagonal matrix with a single 1 at the position
  corresponding to the $x_j$ resp.\ $y_j$ coordinate.
The theorem now follows from the Weyl character formula.
\end{proof}



\section{Even special orthogonal group}
  Let $G=SO_{2r}\left(\C\right)$ be the connected group preserving the quadratic form
$$Q=x_1y_1+\dots+x_ry_r.$$ 
We take the maximal torus $H=
  \left(\C^*\right)^{r}$ consisting of all special orthogonal  diagonal 
matrices $\diag \left(z_1,z_1^{-1},\dots,z_r,z_r^{-1}\right)$.
In  the weight lattice  $
\Hom \left(H,\C^*\right)= \Z^{r}$, we take
  $\lambda=\left(\lambda_1,\dots,\lambda_r\right)$ to correspond to the monomial
  $z^\lambda=\prod z_j^{\lambda_j}$.
The Weyl group  $\mathfrak W$ is the semidirect product of $\mathfrak S_{r}$
and $ Z_2^{r-1}$. It acts by permuting the indices and by performing an even
number of sign changes. Write $\tilde{\mathfrak W}>\mathfrak W$ for the Weyl group in the
full orthogonal group $O_{2r}\left(\C\right)$. It is the semidirect
product of $\mathfrak S_{r}$
and $ Z_2^{r}$. If $\sigma\in\tilde{\mathfrak W}$, we write $[\sigma]$ for its
image in $\mathfrak S_r$.
Write $\rho=\left(r-1,r-2,\dots,1,0\right)$
for  the
half-sum of the positive  roots. Write $\epsilon=\left(1,1,\dots,1,1\right)$  so that
$e=\epsilon+\rho=\left(r,r-1,\dots,2,1\right)$ is 
regular for $\tilde{\mathfrak W}$.
The Weyl character formula is valid as stated in the introduction.
Both sides are $\mathfrak W$--antisymmetric characters of  $H$.

For $\ell\in\Z^{r}$ and $g\in G$, define
  $$g^\ell=\bigotimes_{i=1}^rg^{\ell_i}\in M_{2r}\left(\C\right)^{\otimes r}.$$
We have \begin{align*}
2\sum_{\sigma\in\mathfrak W} \left(-1\right)^\sigma
  g^{\sigma\ell}&=\sum_{\sigma\in\tilde{\mathfrak W}}
  \left(\left(-1\right)^{[\sigma]}+\left(-1\right)
^\sigma\right)
  g^{\sigma\ell}=
\\
&=\bigwedge_{i=1}^{r}
  \left(g^{
\ell_i}+g^{-\ell_i}\right)+\bigwedge_{i=1}^{r}
  \left(g^{
\ell_i}-g^{-\ell_i}\right)\in M_{2r}\left(\C\right)^{\wedge r}.
\end{align*}
Note that the second term is zero if any $\ell_i$ is zero.

\newcommand{\Pf}{\mathrm {Pf}}

\begin{theorem}
Let
  $\lambda
=\left(\lambda_1,\dots,\lambda_r\right)
\in\Z^{r}$ with $\lambda_1\ge\dots\ge\la_{r-1}\ge |\lambda_r|$. Set
  $\bar\la=\left(\lambda_1,\dots,\la_{r-1},-\lambda_r\right)$. Let
 $g\in SO_{2r}\left(\C\right)$. Then

  $$2\sum_{\sigma\in\tilde{\mathfrak W}} \left(-1\right)^{[\sigma]}
  g^{\sigma\ell}=\left(\chi_\lambda+\chi_{\bar\la}\right) \left(g\right)\cdot
 \sum_{\sigma\in\tilde {\mathfrak W}} \left(-1\right)^{[\sigma]}  g^{\sigma\rho};$$
equivalently,
\begin{align*}
2\bigwedge_{i=1}^r\left(g^{\ell_i}+ g^{-\ell_i }\right)=
\left(\chi_\lambda+\chi_{\bar\la}\right) \left(g\right)
\cdot
\bigwedge_{i=1}^{r} 
\left(g^{r-i}+g^{- \left(r-i\right)}\right).
\end{align*}
Also,
  $$\left(\chi_\epsilon-\chi_{\bar\epsilon}\right)\left(g\right)\cdot\sum_{\sigma\in\tilde{\mathfrak W}} \left(-1\right)^{\sigma}
  g^{\sigma\ell}= \left(\chi_\lambda-\chi_{\bar\la}\right) \left(g\right)\cdot
 \sum_{\sigma\in\tilde {\mathfrak W}} \left(-1\right)^{\sigma}  g^{\sigma e};$$
equivalently,
\begin{align*}{\sqrt{-1}}^r\Pf \left(g-g^{-1}\right)
\cdot
\bigwedge_{i=1}^r\left(g^{\ell_i}- g^{-\ell_i }\right)=
\left(\chi_\lambda-\chi_{\bar\la}\right) \left(g\right)
\cdot
\bigwedge_{i=1}^{r} 
\left(g^{r+1-i}-g^{- \left(r+1-i\right)}\right).
\end{align*}  
Throughout, $\ell=\la+\rho$, i.e., $\ell_i=\la_i+r-i$.
\end{theorem}
Note that 
$$\left(\chi_\epsilon-\chi_{\bar\epsilon}\right)\left(g\right)={\sqrt{-1}}^r\Pf \left(g-g^{-1}\right),$$
where the sign of the Pfaffian of {\hbox{$g-g^{-1}\in\mathfrak{so}_{2r}(\C)$}} is specified
by declaring the ordered $Q$--orthonormal bases of the standard vector
representation $\C^{2r}$ with determinant $\left(2\sqrt{-1}\right)^r$
to be of positive orientation.

\begin{proof}
The set of diagonalizable matrices is dense in $G$, so we may
assume that $g$ is such.  The statement of the
theorem is invariant under
conjugation, so we may 
assume that $g=\diag \left(z_1,z_1^{-1},\dots,z_r,z_r^{-1}\right)\in H$. Then 
$$\sum_{\sigma\in\tilde{\mathfrak W}} \left(-1\right)^{[\sigma]}
  g^{\sigma\ell}=\bigwedge_{i=1}^{r}
  \left(g^{\ell_i}+g^{-\ell_i}\right)=
\left|z_j^{\ell_i}+z_j^{-\ell_i}\right|\cdot
\bigwedge_{j=1}^r
  \left(e_{jj}+f_{jj}\right)
$$ 
and 
$$\sum_{\sigma\in\tilde{\mathfrak W}} \left(-1\right)^\sigma
  g^{\sigma\ell}=\bigwedge_{i=1}^{r}
  \left(g^{\ell_i}-g^{-\ell_i}\right)=
\left|z_j^{\ell_i}-z_j^{-\ell_i}\right|\cdot
\bigwedge_{j=1}^r
  \left(e_{jj}-f_{jj}\right),$$
where $e_{jj}$ resp.\ $f_{jj}$
is the diagonal matrix with a single 1 at the position
  corresponding to the $x_j$ resp.\ $y_j$ coordinate.
The theorem now follows from the Weyl character formula.
\end{proof}




\section*{References}

[C] J.\ L.\ Cisneros-Molina, An invariant of $2\times 2$ matrices, Electronic
Journal of Linear Algebra 13 (2005), 146--152.

\bigskip\noindent
[FH] W.\ Fulton and J.\ Harris, Representation theory, GTM, Springer, New
York, 1991.

\end{document}